\documentclass[11pt]{amsart}
\usepackage{amssymb}

\newtheorem{thm}{Theorem}[section]
\newtheorem{cor}[thm]{Corollary}
\newtheorem{prop}[thm]{Proposition}
\newtheorem{lemma}[thm]{Lemma}

\theoremstyle{remark}
\newtheorem{remark}[thm]{Remark}
\newtheorem{example}[thm]{Example}
\newtheorem{question}[thm]{Question}

\theoremstyle{definition}
\newtheorem{defn}[thm]{Definition}

\numberwithin{equation}{section}

\newcommand{\bbA}{{\mathbb A}}
\newcommand{\bbC}{{\mathbb C}}

\newcommand{\bbZ}{{\mathbb Z}}

\newcommand{\bbP}{{\mathbb P}}
\newcommand{\bbN}{{\mathbb N}}
\newcommand{\bbF}{{\mathbb F}}

\newcommand{\bbQ}{{\mathbb Q}}

\newcommand{\A}{{\mathbb A}}

\newcommand{\N}{{\mathbb N}}

\newcommand{\G}{{\mathbb G}}

\newcommand{\ze}[1]{\zeta _{#1\,\sigma}(t)}

\newcommand{\cF}{{\mathcal F}}

\newcommand{\cO}{{\mathcal O}}

\newcommand{\cL}{{\mathcal L}}

\newcommand{\cV}{{\mathcal V}}

\newcommand{\V}{{\mathcal V}}

\newcommand{\Ind}{\operatorname{Ind}}
\newcommand{\Spec}{\operatorname{Spec}\,}
\newcommand{\Hilb}{{\rm Hilb}}
\newcommand{\Sym}{{\rm Sym}}
\newcommand{\Stab}{{\rm Stab}}
\newcommand{\Sign}{{\rm Sign}}
\newcommand{\Vect}{{\rm Vect}}
\newcommand{\groth}[1]{K_0[\V_{#1}]}
\newcommand{\grsig}[1]{K_0[\V_{#1}]_\sigma}

\newenvironment{pf}{\smallskip\noindent{\bf Proof.}\ }{\qed\smallskip}

\title{Rationality criteria for motivic zeta-functions}
\author{Michael Larsen}
\address{Department of Mathematics, Indiana University,
Bloomington, IN 47405, USA}
\email{larsen@math.indiana.edu}
\author{Valery A.~Lunts}
\address{Department of Mathematics, Indiana University,
Bloomington, IN 47405, USA}
\email{vlunts@indiana.edu}

\thanks{The first named author was partially supported by NSF grant DMS-0100537.  
The second named author was partially supported by NSA grant MDA904-01-1-0020  
and CRDF grant RM1-2405-MO-02}

\begin{document}
\maketitle

\section{Introduction}

Let $X$ be a variety over a finite field $\bbF$ and $\Sym^n X = X^n/S_n$, the variety
of effective zero-cycles on $X$ of degree $n$ (where by convention,
$S^0 X = \Spec\bbF$).  A celebrated theorem of
B.~Dwork asserts that the zeta-function
$$Z_X(t) = \sum_{i=0}^\infty |\Sym^n(X)(\bbF)| t^n$$
is a rational function in $t$.  Kapranov asked \cite{Ka} whether this 
rationality lifts to the Grothendieck ring of varieties over $\bbF$.
Explicitly, let
$\groth{\bbF}$ denote the ring of $\bbZ$-combinations of isomorphism
classes of $\bbF$-varieties modulo the cutting-and-pasting relation
$$[X] = [Y] + [X\setminus Y]$$
for closed $\bbF$-subvarieties $Y\subset X$.   Is the {\bf motivic zeta-function}
$$\zeta_X(t) = \sum_{i=0}^\infty [\Sym^n(X)] t^n$$
always rational as a power series in $\groth{\bbF}$?   More generally,
is this true for varieties over a general field $K$, for example, $K=\bbC$?
If so, this would give a fundamentally new proof of Dwork's
theorem, one that does not depend on Frobenius at all.  Kapranov observed that it
{\it is} true when $X$ is a curve.  (See Theorem~\ref{onedim} below for a statement and proof
over $\bbC$.)  

In [LaLu], we proved that $\zeta_X(t)$ is {\it not} rational when $X$ is a complex projective
non-singular surface with geometric genus $\ge 2$.  In this paper we return to
the problem and give a simple necessary and sufficient condition on complex surfaces
for $\zeta_X(t)$ to be rational:

\begin{thm}
A complex surface $X$ has rational motivic zeta function if and only if it has 
Kodaira dimension $-\infty$.
\end{thm}

The notion of rationality needs some explanation (see \S2 below).
Theorem~\ref{positive} asserts rationality in the strong sense when
Kodaira dimension is $-\infty$, and Theorem~\ref{negative} denies rationality in
the weak sense when Kodaira dimension is $\ge 0$.
The methods for proving rationality and irrationality are entirely different.  The
fact that they meet in the middle to give a necessary and sufficient
condition appears as a minor miracle,
made possible by the classification of surfaces.
By contrast, in dimension $\ge 3$, the rationality problem seems wide open.

There are two main ideas in this paper.
The first is to probe rationality by means of new {\bf motivic measures},
i.e. field-valued points of $\Spec\groth \bbC$.
In [LaLu], we showed that
any multiplicative function on the set of stable birational equivalence
classes of non-singular complex projective varieties defines a motivic measure.  
In particular, we constructed a measure $\mu_1$, with values in the 
fraction field of the group ring of the
multiplicative group of integer power series $(1+s\bbZ[[s]])^\times$,
characterized by the formula
$$\mu_1([Z]) = \Bigl[\sum_{i=0}^\infty \dim\Gamma(Z,\Omega_Z^i) s^i\Bigr]$$
for $Z$ a (connected) non-singular projective variety.
In the current paper, we introduce a set of measures $\mu_n$ indexed by positive
integers $n$, characterized by
\begin{equation}
\label{measures}
\mu_n([Z]) = \Bigl[\sum_{i=0}^\infty \dim H^0\bigl(Z,\Psi^n\Omega_Z^i\bigr) s^i\Bigr],
\end{equation}
where $n$ denotes the $n$th Adams operation. 
The point of this generalization is that it allows us to work with higher plurigenera 
in much the way that we worked with geometric genus in [LaLu].  If for some $n$,
$\mu_n(\zeta_X(t))$ is irrational, then of course $\zeta_X(t)$ itself is irrational.
For singular $Z$, the right hand side of (\ref{measures}) does not make sense, 
and indeed the left hand side is generally not of the form $[P]$ for 
any power series $P\in (1+s\bbZ[[s]])^\times$.  However, in the special case
$Z = \Sym^n X$, $\dim X= 2$,
it turns out that $\mu_n([Z])= \mu_n([\Hilb^n X])$, where $\Hilb^n X$ is the
Hilbert scheme of points on $X$.  (The fact that $\Hilb^n X$ is non-singular
is another way in which 
dimension $2$ is special.)   Our task, therefore, is to prove that if the $m$th plurigenus
of $X$ is positive, then
$$\sum_{i=0}^\infty\Bigl[\sum_{j=0}^\infty 
\dim H^0\bigl(\Hilb^i X,\Psi^m\Omega_{\Hilb^i X}^j\bigr) s^j\Bigr] t^i$$
is irrational.

The other main idea in this paper is to make systematic use of $\lambda$-ring
ideas and techniques (see \S4 for a brief review of this theory).  
We have already noted the appearance of Adams operations.   This is somewhat delicate:
we cannot work in the usual $K(X)$ since the global section functor is
not well-defined there.  Instead, we need to prove that the ring of virtual vector
bundles modulo {\it split} exact sequences is a special $\lambda$-ring (see \S5).
Also significant is the idea that the
motivic zeta-function (or better, its inverse) should be regarded
as the universal $\lambda$-homomorphism. 
This is true only formally since $\groth{\bbC}$ is not a special $\lambda$-ring.
It is nonetheless suggestive, since the universal $\lambda$-homomorphism
sends every virtually finite element to a rational power series.  For a special 
$\lambda$-ring the virtually finite elements form a $\lambda$-subring which
in many interesting cases is the whole ring.  This interpretation of the zeta-function
suggests on the one hand that we should seek motivic measures which are $\lambda$-homomorphisms
and on the other that it may be natural to replace the Grothendieck ring of varieties
by its {\it specialization}, i.e., its maximal quotient which is a special $\lambda$-ring
(see \S8).

E.~Looijenga has called our attention to certain formal analogies between the problem of
rationality of zeta-functions of complex surfaces $X$, and Severi's conjecture,
disproved by D.~Mumford \cite{Mu}, which predicted
that the group of $0$-cycles modulo rational equivalence on $X$
would be finite-dimensional.  Each statement can be understood as bounding the
complexity of symmetric powers of $X$.  In each case, global holomorphic $2$-forms 
provide the obstruction to this boundedness.  However, there is an essential difference.
A conjecture of S.~Bloch \cite{Bl} asserts that when $H^0(X,\Omega^2_X) = 0$,
the group of degree zero $0$-cycles should be isomorphic to the Albanese variety of $X$.
This has been settled by Bloch, A.~Kas, and D.~Lieberman
\cite{BKL} for Kodaira dimension $<2$, so there are examples of surfaces
in which the two notions of boundedness diverge.  In particular, the use of plurigenera
in this paper does not seem to have a counterpart in the world of ${\rm CH}^2$.

Throughout this paper, a {\bf variety} will be a reduced separated scheme 
of finite type over $\bbC$.
The class $[X]\in\groth{\bbC}$ of a variety $X$ will sometimes be written without brackets when it 
seems unlikely to lead to any confusion.

\section{Rationality criteria for power series}

In this paper we will be concerned with the rationality of power series
with coefficients
in a commutative ring $A$.  It is not entirely clear how such rationality
should be defined when $A$ is not an integral domain (and we know \cite{Po}
that the Grothendieck ring of varieties has zero-divisors).
We consider several possible definitions:

\begin{defn}
\label{globrat}
A power series $f(t)\in A[[t]]$ is {\bf globally rational} if and
only if there exist
polynomials $g(t), h(t)\in A[t]$ such that $f(t)$ is the unique solution
of $g(t)x = h(t)$.
\end{defn}

\begin{defn}
\label{detrat}
A power series $f(t)=\sum_i a_i t^i\in A[[t]]$ is
{\bf determinantally rational}
if and only if there exist integers $m$ and $n$ such that
$$\det\left(\begin{array}{cccc}
a_i&a_{i+1}&\cdots&a_{i+m}\\
a_{i+1}&a_{i+2}&\cdots&a_{i+m+1}\\
\vdots&\vdots&\ddots&\vdots\\
a_{i+m}&a_{i+m+1}&\cdots&a_{i+2m}
\end{array}\right)=0$$
for all $i > n$.
\end{defn}

It is classical \cite{E.Bo} that Definition~\ref{globrat} is equivalent to
Definition~\ref{detrat} when $A$ is a field.  This suggests a third
possible definition:

\begin{defn}
A power series $f\in A[[t]]$ is {\bf pointwise rational} if and
only if
for all homomorphisms $\Phi$ from $A$ to a field, $\Phi(f)$ is rational by
either of the two previous definitions.
\end{defn}

These definitions are related by the following proposition:

\begin{prop}
Any globally rational power series is determinantally rational, and any
determinantally
rational power series is pointwise rational.  Neither converse holds for a
general
coefficient ring $A$.  All three conditions are equivalent when $A$ is
an integral domain.

\end{prop}

\begin{pf}
Suppose $g(t)=\sum_i^k b_i t^i,h(t)=\sum_i^\ell c_i t^i\in A[[t]]$
and $f(t)=\sum_i a_i t^i$ is the unique solution to
$g(t) x = h(t)$.   Since $f(t)+a$ is not a solution when $a\neq 0$
is a constant, the annihilator of the ideal $(b_0,b_1,\ldots,b_k)$
of coefficients of $g$ is $(0)$.  Setting $m=k$, $n=\ell$ in Definition~\ref{detrat},
$$\left(\begin{array}{cccc}
a_i&a_{i+1}&\cdots&a_{i+m}\\
a_{i+1}&a_{i+2}&\cdots&a_{i+m+1}\\
\vdots&\vdots&\ddots&\vdots\\
a_{i+m}&a_{i+m+1}&\cdots&a_{i+2m}
 \end{array}\right)\left(\begin{array}{c}
b_k\\
b_{k-1}\\
\vdots\\
b_0
 \end{array}\right)=\left(
\begin{array}{c}
0\\
0\\
\vdots \\
0
\end{array}
\right)$$
for $i>\ell-k$.
Left multiplying both sides by the matrix of cofactors, we conclude that
the determinant of the above matrix annihilates the $b$-column vector,
which means that
it is $0$.

To see that the converse does not hold, consider $A=\bbZ [x]/(x^2)$.
As $A$
is countable, the set of globally rational power series over $A$ is
countable.
However, any power series of the form $f(t) = xg(t)$ satisfies
Definition~\ref{detrat} when $m=1$.

Suppose $f(t)=\sum_i a_i t^i$ satisfies Definition~\ref{detrat}.  The same is then
true for
$\Phi(f)$, which is a power series over a field.  For fields, however, the
determinantal
condition implies rationality.

To see that the converse does not hold, consider
$$A=\bbZ [x_1,x_2,\ldots]/(x_1^2,x_2^2,\ldots)$$
in infinitely many variables $x_i$.
Every homomorphism from $A$ to a field factors through the augmentation.
Therefore,
$f(t)=\sum_{i=1}^\infty x_i t^i$ is pointwise rational.  However, the
determinant
$$\det\left(\begin{array}{cccc}
a_i&a_{i+1}&\cdots&a_{i+m}\\
a_{i+1}&a_{i+2}&\cdots&a_{i+m+1}\\
\vdots&\vdots&\ddots&\vdots\\
a_{i+m}&a_{i+m+1}&\cdots&a_{i+2m}
\end{array}\right)
$$
is never zero since it has a non-trivial $x_i x_{i+2} x_{i+4} \cdots
x_{i+2m}$
coefficient.

Finally, when $A$ is an integral domain, we let $\Phi$ denote the
inclusion map from $A$ to its fraction field $F$.  If $f(t)$ is pointwise 
rational, it must be rational as a power series in $F[[t]]$.  Therefore,
$g(t) f(t) = h(t)$ for polynomials $g(t),h(t)\in F[t]$, $g(t)\neq 0$.
Clearing denominators, we may assume  $g(t),h(t)\in A[t]$.  As $A[t]$
is an integral domain, $f(t)$ is globally rational.
\end{pf}

The following lemma will be useful in \S3:

\begin{lemma}
\label{invrat}
If $f(t)\in 1+tA[[t]]$ is globally (resp. pointwise) rational,
the same is true of $f(t)^{-1}$.
\end{lemma}

\begin{pf}
If $f(t)$ is the unique solution of $g(t) x = h(t)$, then $g(t)$ is not a
zero-divisor in
$A[[t]]$, and $f(t)$ is invertible, so again not a zero-divisor.
Therefore $h(t)=f(t)g(t)$
is not a zero divisor, and $f(t)^{-1}$ is the unique solution of $h(t) x =
g(t)$.
The pointwise case is trivial.
\end{pf}

In \cite{Ka}, the rationality of motivic zeta functions is discussed in
pointwise terms.
To give the strongest possible results, we
generally prove rationality globally and irrationality pointwise.  
To do the latter, we
make free use of the determinantal formulation of rationality for power
series over fields.

\medskip

 Consider a free  abelian group $G$ and its group ring $\bbZ [G]$. This
ring is isomorphic to the ring of Laurent polynomials, hence is a domain. Let $F$ be
its field of fractions. Denote by $\Theta$ the collection of power series
$\sum g_it^i\in F[[t]]$, where $g_i\in G$ or $g_i=0$.  The following proposition
gives a rationality criterion for power series in $\Theta$.

\medskip

\begin{prop} 
\label{monomial}
Let $f(t)=\sum g_it^i\in \Theta$ be a  power series in $F[[t]]$. Then
$f(t)$ is rational if and only if
there exists $n\geq 1$, $i_0\in\bbN$ and a sequence $\{h_i\in G\}$ periodic with period $n$ 
such that
for $i> i_0$ we have $g_{i+n}=h_ig_i$.
\end{prop}

\begin{pf} 
If there exist $n$, $i_0$, and $h_i$ as above, then
$$f(t) = \sum_{i=0}^{i_0} g_i t^i + \sum_{i = 1}^n \frac{g_{i_0+i}}{1-h_i t^n}$$
is a rational function with denominator of degree $m\le n^2$.

For the converse, by assumption there exists a polynomial
$$q(t)=
a_mt^m+a_{m-1}t^{m-1}+\cdots+a_0 \in F[t]$$
such that
$q(t)f(t)$ is a polynomial. 
Clearing denominators, we may assume $a_i\in\bbZ[G]$ for all $i$.
Let $C\subset G$ be the collection $\alpha \beta ^{-1}$ where
$\alpha$, $\beta $ are group elements which
have non-zero coefficient in some $a_j$, $0\le j\le m$. Then for $i> i_0$
\begin{equation}
\label{recursion}
-a_0g_{i+m}=a_1g_{i+m-1}+\cdots+a_mg_i
\end{equation}
and in particular
\begin{equation}
\label{inclusion}
g_{i+m}\in g_iC \cup \cdots\cup g_{i+m-1}C \cup \{0\}
\end{equation}

\medskip

\begin{defn} A power series $p(x)=\sum\alpha _it^i\in \Theta$ is {\it compact} 
(or {\it $(m,K)$-compact}) if there
exists a finite set $K\subset G$ such that for nonzero $\alpha _i$, $\alpha _j$ such that
$|i-j|\leq m $ we have $\alpha _i\alpha _j^{-1}\in K$. 
\end{defn}

\medskip

\noindent{\it Case 1.} Assume that $f(t)$ is compact. Fix $i > i_0$. Using the compactness
of $f(t)$ we can find $n\geq 1$ and $h\in G$ such that
$$(g_{i+n},g_{i+n+1},\ldots,g_{i+n+m-1})=h(g_{i},g_{i+1},\ldots,g_{i+m-1}).$$
Then (\ref{recursion}) implies
$$(g_{j+n},g_{j+n+1},\ldots,g_{j+n+m-1})=h(g_{j},g_{j+1},\ldots,g_{j+m-1})$$
for all $j\geq i$, which proves the proposition in this case. Note that we can bound
 $n$ by some function $\phi(m,|K|)$.

\medskip

\noindent{\it Case 2.} This is the general case. It follows from (\ref{inclusion}) that
there exists an integer $r$, $1\leq r\leq m$ and compact
power series $f_1(t),\ldots,f_r(t)$ in
$\Theta $ such that $f(t)=\sum f_j(t)$. We choose the {\it minimal} such $r$ and corresponding
series $f_j(t)=\sum g^j_it^i$ for $j=1,2,\ldots,r$. 

\medskip

\begin{lemma} The power series $f_j(t)$ simultaneously satisfy the recurrence relations
(\ref{recursion}) for coefficients with indices contained in arbitrarily long intervals.
\end{lemma}

\begin{pf} We use the minimality of $r$. If $r=1$ we are in Case 1 and there is
nothing to prove. Assume $r\geq 2$.   The minimality of $r$ implies (and actually is
equivalent to) the following:

For every finite subset $S\subset G$ there exist infinitely many intervals $I\subset \bbN$
of length $m+1$ such that for any nonzero $g_{i_1}^k$, $g_{i_2}^l$ with $i_1,i_2\in I$ and
$k\neq l$ we have $g_{i_1}^k(g_{i_2}^l)^{-1}\notin S$.

In particular, we can find such an interval $I=(i,i+1,\ldots,i+m)$ for the set $S=C$ so that
$i > i_0$. Then the recurrence relation (\ref{recursion}) implies the same relation for
coefficients $g_i^j,\ldots,g_{i+m}^j$ for all $j=1,\ldots,r$. To get a long interval $J$,
say of length $sm$ on which this recurrence holds for all $j$ we proceed as follows:
let $K\subset G$ be finite, such that each $f_j(t)$ is $(m,K)$-compact; take $S=CK^{2s}$ and
let $I=(i,i+1,\ldots,i+m)$ be a corresponding interval such that $i >  i_0$. Then we can take
$J=(i,i+1,\ldots,i+sm)$.
This proves the lemma.
\end{pf}

We now complete the proof of the proposition as follows: 
Choose a sufficiently long interval $J\subset (i_0,\infty)$ as
in the proof of last lemma.
Then repeating the argument in Case 1 for each $j=1,\ldots,r$ we find $n_j\geq 1$ and
$h^\prime_j\in G$ such that
$$(g_{i+n_j}^j,g_{i+n_j+1}^j,\ldots,g_{i+n_j+m-1}^j)
=h^\prime_j(g_{i}^j,g_{i+1}^j,\ldots,g_{i+m-1}^j)$$
as long as $[i+n_j+m]\subset J$. Now take $n=n_1\cdots n_r$ and put
$$h_i=\begin{cases}(h^\prime _j)^{n/n_j},\ \text{if $i\in J$ and $g_i^j\neq 0$},\\
                    1,\ \text{if $i\in J$ and all $g_i^j=0$}.\\
\end{cases}$$
It follows that $g_{i+n}=h_ig_i$ as long as the indices stay in the interval $J$.

The recurrence
(\ref{recursion}) implies that there exists a linear operator $T$ such that 
$$(g_{i+n},g_{i+n+1},\ldots,g_{i+n+m}) = T^n(g_{i},g_{i+1},\ldots,g_{i+m})$$
for all $i\geq i_0$. We claim that the relation
$g_{i+n}=h_ig_i$ holds in fact for {\it all} $i\geq \inf(J)$. Indeed, this follows from the
following simple fact in linear algebra: Let $V$ be a $d$-dimensional vector space
and $0\neq v\in V$; assume that we have two linear operators $A$ and $B$ in $V$, such that
$A^lv=B^lv$ for $l=1,\ldots,d-1$. Then $A^lv=B^lv$ for all $l\geq 0$.
\end{pf}

\section{Rationality theorems}

In this section we describe several classes of varieties which have
rational zeta-functions,
notably rational and ruled surfaces and linear algebraic groups.  These
are rather special
cases, and indeed we suspect that except in dimension 1, rationality is
the exception rather
than the rule.  Some evidence for this point of view is given in the
discussion of
irrationality theorems in \S7.

\begin{lemma}
\label{anytwo}
If $X$ is a variety and $Y$ a closed subvariety with complement $U$,
then if any two of $\zeta_X(t)$, $\zeta_Y(t)$, and $\zeta_U(t)$
are globally (resp. pointwise) rational, then the third is so as well.
\end{lemma}

\begin{pf}
 Stratifying $\Sym^n X$ according to how many points land in $Y$, we
obtain
\begin{equation}
\label{whylambda}
[\Sym^n X] = \sum_{i+j=n}[ \Sym^i Y \times \Sym^j U] = \sum_{i+j=n}
[\Sym^i Y][\Sym^i U].
\end{equation}
It follows that $\zeta_{X}(t) = \zeta_Y(t)\zeta_U(t)$.

By Lemma~\ref{invrat}, not only are the two specified zeta-functions globally (resp.
pointwise)
rational, the same is true of
their reciprocals.  The product of rational zeta-functions is again
rational,
and the lemma follows.
\end{pf}

In particular, the disjoint union of varieties with globally rational
zeta-functions again
has a globally rational zeta-function.  This shows in particular that
zero-dimensional varieties have globally rational zeta-functions.
Also, a stratified variety has such a
zeta-function as long as all of its strata do.

\begin{lemma}
\label{vbclass}
If $X$ is any variety,\ and $E\to X$ is a vector bundle of rank $r$, then
$$[E] = [X][\A^r].$$
\end{lemma}

\begin{pf}
Vector bundles are locally trivial in the Zariski topology, so there
exists a dense open
subset $U\subset E$ over which $E$ restricts to a trivial bundle.  If $Y$
denotes the complement
of $U$ and $E_U$ and $E_Y$ the pull-back of $E$ to $U$ and $Y$
respectively, then
applying Lemma~\ref{anytwo} to $E_Y\subset E$ and $Y\subset X$, it suffices to
prove the proposition
for $Y$ and $U$, and for $U$, $E_U = U\times\A^r$.  The lemma follows by
Noetherian induction.
\end{pf}

The following proposition is due to B.~Totaro \cite{Go}:

\begin{prop}
If $X$ is any variety, $E\to X$ is a vector bundle of rank $r$,
and $n$ is any positive integer, then $[\Sym^n E] = [\Sym^n X]
[\A^1]^{rn}$.
Equivalently, $\zeta_E(t) = \zeta_X(\A^r t)$.
\end{prop}

\begin{cor}
\label{vbrat}
If $X$ is a variety such that $\zeta_X(t)$ is globally (resp. pointwise)
rational,
and $E\to X$ is a vector bundle,
then $\zeta_E(t)$ is globally (resp. pointwise) rational.
\end{cor}

\begin{cor}
\label{cells}
For all non-negatve integers $n$, $\A^n$ is rational.
\end{cor}

\begin{cor}
\label{proj}
If $X$ is a variety and $P\to X$ is a projective space bundle of rank $r$
which is locally trivial in the Zariski topology, then
$$\zeta_P(t) = \zeta_X(t)\zeta_X(\A^1 t)\cdots\zeta_X(\A^n t).$$
\end{cor}

In particular, $\zeta_{{\mathbb P}^r}(t)$ is rational for all $r$.

\begin{pf}
By Noetherian induction, it suffices to consider the case $P =
X\times{\mathbb P}^r$.
The corollary then follows immediately from the stratification of ${\mathbb P}^r$
with strata $\A^i$, $0\le i\le r$.\end{pf}

%
%
%
%

We come now to the main positive result in the subject.

\begin{thm}(Kapranov) If $X$ is any one-dimensional variety, then $\zeta_X(t)$ is
globally rational.
\label{onedim}
\end{thm}

\begin{pf}
Let $X$ be a $1$-dimensional variety.  The singular locus $Y$ is either
empty or zero-dimensional.
In the latter case, the rationality question reduces to the case of
$X\setminus Y$, so without loss of generality we may assume $X$ is
non-singular.  If $X$ has
more than one component, it suffices to prove that each one has a rational
zeta-function,
so without loss of generality we may assume $X$ is connected and therefore
irreducible.
Let $\bar X$ denote the unique projective non-singular curve containing
$X$ as an open subvariety.
As $\bar X\setminus X$ is empty or zero-dimensional, without loss of
generality we may
assume that $X$ is projective and non-singular.

Let $g$ be the genus of $X$.  Let $x_0$ denote a base point of $X$ and $J$
denote
the Jacobian variety ${\rm Jac}^0(X)$.
For non-negative $n\ge 2g-1$,  the morphism
$X^n\to J$ mapping
$$(x_1,\ldots,x_n)\mapsto -n x_0 + \sum_{i=1}^n x_i$$
factors through $\Sym ^nX$ realizing it as
a projective
space bundle of rank $n-g$ over $J$.  The closed immersion $X^n\to
X^{n+1}$
sending $(x_1,\ldots,x_n)$ to $(x_0,\ldots,x_n)$ induces a closed
$J$-immersion $\Sym^n X\to \Sym^{n+1}(X)$; the complement is a vector
bundle of
rank $n+1-g$ over $J$.  By Lemma~\ref{vbclass},
$$[\Sym^{n+1} X] - [\Sym^n X] = [X][\A^1]^{n+1-g}.$$
This implies that
\begin{equation}\label{curveNumer} \zeta_X(t)(1-t)(1-\A^1 t)
\end{equation}
is a polynomial of degree $\leq 2g$.
\end{pf}

\begin{cor}
If $X$ is an algebraic surface, the global (resp. pointwise)
rationality of $\zeta_X(t)$ depends only
on the birational equivalence class of $X$.
\end{cor}

\begin{pf}
If $X_1$ and $X_2$ are two such surfaces and $U$ is a surface which is
isomorphic
to dense open subvarieties of each, then setting $Y_i = X_i\setminus U$,
$$\zeta_{X_i}(t) = \zeta_U(t)\zeta_{Y_i}(t).$$
By Theorem~\ref{onedim}, $\zeta_{Y_i}(t)$ is globally (therefore also pointwise)
rational, so by Lemma~\ref{anytwo}, each $\zeta_{X_i}(t)$ is rational if and only
if $\zeta_U(t)$
is so.\end{pf}

\begin{thm}
\label{positive}
If $X$ is a surface with Kodaira dimension $-\infty$, then $\zeta_X(t)$ is
globally rational.
\end{thm}

\begin{pf}
There are two cases: rational surfaces and birationally ruled surfaces.
In each case, we may choose
any variety in the given birational equivalence class.  For rational
surfaces, we use
$\A^2$, which has a globally rational zeta-function by Corollary~\ref{cells}.  A ruled
surfaces is a
projective line bundle over a curve.  By the Tsen-Lang theorem, such a
${\mathbb P}^1$-bundle is
Zariski-locally trivial, and by Theorem~\ref{onedim} and Corollary~\ref{proj}, the motivic
zeta-function of a
ruled surface is globally rational.
\end{pf}

\begin{thm}
If $X$ is a linear algebraic group, its zeta-function is globally
rational.
\end{thm}

\begin{pf}
The components of an algebraic group $G$ are isomorphic to one another as
varieties, so
$$\zeta_G(t) = \zeta_{G^\circ}(t)^{[G:G^\circ]},$$
where $G^\circ$ denotes the identity component of $G$.  We therefore
assume that $G$ is
connected.  Let $U$ denote the unipotent radical $G$ and $H=G/U$.
Then $U$ is isomorphic to a
subgroup of the unitriangular matrices $n\times n$ matrices  \cite{A.Bo}~4.8.  The logarithm
map therefore defines an isomorphism of varieties between $U$ and its Lie algebra, which 
means that $U$ is isomorphic to $\A^{\dim U}$.  The existence of Levi decompositions gives a 
an isomorphism $G\tilde\to U\times H$, so it suffices to prove
that $\zeta_H(t)$ is globally rational.

Let $B$ denote a Borel subgroup of $H$, $T$ a maximal torus of $B$,
and $V$ the unipotent radical of $B$.  Decomposition into Schubert cells
gives a stratification of the flag variety $H/B$ in which every stratum is
isomorphic
to $\A^k$ for some $k$ with pre-image $V w B\subset H$ isomorphic 
to $\A^k\times B \cong \A^{k+\dim V}\times T$ \cite{A.Bo}~14.12.
It therefore suffices to prove that $\zeta_T(t)$ is rational.

We use induction on $d = \dim T$.  The rationality is trivial for $d=0$.  Now,
$$\G_m^d \times \A^1 \setminus \G_m^d = \G_m^{d+1},$$
so by Corollary~\ref{vbrat} the rationality of $\zeta_{\G_m^d}(t)$ implies that of
$\zeta_{\G_m^{d+1}}(t)$.

\end{pf}

\section{Lambda rings}
In this section, we develop some basic definitions and facts connected
with
the notion of $\lambda$-ring.  A good reference for this material is the
first section
of \cite{AtiTa}

\begin{defn}
A $\lambda$-{\bf structure} on a commutative ring $A$ is an infinite
sequence
$\lambda^0,\lambda^1,\lambda ^2,\ldots,$ of maps $A\to A$ such that
$$\begin{array}{rcl}
\lambda^0(x)&=&1\\
           \lambda^1(x)&=&x\\
           \lambda^n(x+y)&=&\sum_{i+j=n}\lambda^i(x)\lambda^j(y).
\end{array}
$$
A $\lambda$-{\bf ring} is a commutative ring endowed with a
$\lambda$-structure.
We call a ring homomorphism between
$\lambda$-rings which
commutes with $\lambda$-operations a $\lambda$-{\bf homomorphism}.
\end{defn}

The prototype of a $\lambda$-ring is the Grothendieck group of finite
dimensional
vector spaces over a field.  Of course, explicitly the ring here is $\bbZ$.
More generally, one can look at Grothendieck groups of finitely generated
projective modules over a ring
or finite rank vector bundles over variety or over a topological space.
In each case, $\lambda^i$
should be regarded as the $i$th exterior power operation.  For example, in
the first case,
$\lambda^i(n)=\binom{n}{i}$.

\begin{example} 
\label{vvs}
Consider the polynomial ring $\bbZ [s]$. If we identify $\bbZ$ with
the ring of virtual finite dimensional vector spaces, then $\bbZ [s]$ is identified
with the ring of isomorphism classes of
$\N$-graded virtual finite dimensional vector spaces. 
Using this
identification define the $\lambda $-structure on $\bbZ [s]$ as follows:
$$\lambda ^i(Vs^p)=\begin{cases}
\Sym ^iVs^{ip},\hbox{ if $p$ is even,}\\
\Lambda ^iVs^{ip},\hbox{ if $p$ is odd.}\\
\end{cases}
$$
Here $\Lambda$ and $\Sym$ have their usual meanings for ordinary virtual vector
spaces:
\begin{align*}
\Lambda^k (V-W) &= \sum_{i+j=k} (-1)^j\Lambda^i V\otimes \Sym^j W;\\
   \Sym^k (V-W) &=  \sum_{i+j=k} (-1)^j\Sym^i V\otimes \Lambda^j W.\\
\end{align*}
%
%
%
Note that the sign conventions for $\lambda^i$ 
are those one expects for {\it symmetric} powers
of virtual super vector spaces. 
\end{example}

\begin{example}
\label{monoidring}
Let $A$ be a $\lambda$-ring and $M\subset A$ a multiplicative submonoid
closed under the
$\lambda$-operations. Consider the corresponding monoid ring $\bbZ [M]$. (Note that
it is {\it not} a subring of $A$).
Then $\bbZ[M]$ has a natural $\lambda$-structure
given by
$$\lambda^i([m])=[\lambda^i m].$$
\end{example}

\begin{example} 
\label{target}
Let us combine the two previous examples. Let $M\subset \bbZ [s]$ be
the multiplicative monoid, which consists of polynomials with constant term 1.
Then $\bbZ [M]$ is a $\lambda $-ring. This example will be important to us since
our motivic measures will take their values in $\bbZ[M]$. We will also need the following
lemma.
\end{example}

\medskip

\begin{lemma} 
\label{monoidm}
The monoid $M$ is a free commutative monoid. The ring $\bbZ[M]$ is a
polynomial ring, hence an integral domain.
\end{lemma}

\begin{pf} The ring $\bbZ[s]$ is factorial and any element of $M$ is a unique product
of elements of $M$, which are prime in $\bbZ[s]$ (the only unit in $M$ is 1). Thus
$M$ is isomorphic to the monoid $\oplus \N$, where the summation is over all prime elements of
$\bbZ[s]$ which are in $M$. Hence $\bbZ[M]$ is a polynomial ring, so it is an integral domain.
\end{pf}

\begin{example} 
\label{source}
By (\ref{whylambda}), the Grothendieck ring of varieties $\groth{K}$ has a natural $\lambda $-ring structure for any field $K$. Indeed, put $\lambda ^i[X]:=[\Sym ^iX]$. Then
$$\lambda _t([X])=\sum_{i=0}^{\infty}\lambda ^i[X]t^i=\zeta _X(t).$$
\end{example}

For any commutative ring $A$ there is a natural $\lambda$-ring structure
on
the set $1+t A[[t]]$.  The operation of addition in this ring is given
by the multiplication
of power series.  Multiplication ($\cdot$) and the $\lambda$-operations
($\Lambda^i$) are
given by universal polynomials which are uniquely characterized by the
identities
\begin{equation}
\label{prodform}
\Bigl\{\prod_{i=1}^m (1+a_i t)\Bigr\}\cdot\Bigl\{\prod_{j=1}^n
(1+b_jt)\Bigr\}
=\prod_{i=1}^m\prod_{j=1}^n (1+a_i b_jt)
\end{equation}
and
$$\Lambda^k\prod_{i=1}^n (1+a_i t)
= \prod_{
S\subset\{1,\ldots,n\},\;|S|=k
}
(1+t\prod_{j\in S} a_j).
$$
Equation (\ref{prodform}) implies that the polynomial expressing the $t^p$ coefficient of
$\sum_i x_i t^i \cdot \sum_j y_j t^j$ in terms of $x_i$ and $y_j$ lies in the ideal
$$(x_{m+1},x_{m+2},\ldots,y_{n+1},y_{n+2},\ldots)$$
whenever $p>mn$.  This implies that the product of two polynomials in $1 + t A[[t]]$
is again a polynomial, regardless of whether the polynomials split into linear factors.
As $\cdot$ distributes over the usual multiplication of power series, $f(t)\cdot g(t)$
is a ratio of polynomials in $1 + t A[[t]]$ if $f(t)$ and $g(t)$ are.

If $A$ has a $\lambda $-structure then the map
$$\lambda_t\colon A\to 1+t A[[t]],\ \lambda_t(a)=\sum_{i=0}^\infty
\lambda^i a t^i$$
is an additive group homomorphism.

\begin{defn}
\label{fdim}
An element $a\in A$ is {\bf finite dimensional} if $\lambda_t(a)$ is a
polynomial,
and the {\bf dimension} of $a$ is the degree of this polynomial.
A difference of finite dimensional elements is {\bf virtually finite}.
A $\lambda$-ring is {\bf finite dimensional} if all of its elements are
virtually finite.
\end{defn}

\begin{defn}
We say $A$ is a {\bf special} $\lambda$-ring if the homomorphism of additive groups
$\lambda_t$ is a
$\lambda$-homomorphism.
In this case, we will call $\lambda$ the universal
$\lambda$-homomorphism.
\end{defn}

For example, the Grothendieck group of vector spaces (or projective
modules or vector bundles)
is a special $\lambda$-ring.  Moreover,
$1+t A[[t]]$ is always special, regardless of whether or not $A$ is
so.
On the other hand, Examples \ref{vvs}, \ref{target}, and \ref{source} are not special,
and neither is Example~\ref{monoidring} in general, even if
$A$ happens to be so.  Special $\lambda$-rings are characterized by
identities of the form
\begin{equation}\label{ringhom}
\lambda^n(xy) = P_n(x,\lambda^2 x,\ldots,\lambda^n x,y,\ldots,\lambda^n
y)\end{equation}
and
\begin{equation}\label{lambdahom}
\lambda^m\lambda^n x = Q_{m,n}(x,\lambda^2 x,\ldots,\lambda^{mn}
x)\end{equation}
for certain universal polynomials $P_n$, $Q_{m,n}$.  Equation
~(\ref{ringhom})
guarantees
that $\lambda_t$ is a ring homomorphism and (\ref{lambdahom}) guarantees it
respects $\lambda$-structures.

\begin{lemma}
For any $\lambda$-ring $A$ there exists a  universal pair consisting of a
special $\lambda$-ring
$B$ and a $\lambda$-homomorphism $A\to B$ such that every
$\lambda$-homomorphism from $A$
to a special $\lambda$-ring $C$ factors through $B$.
\end{lemma}

\begin{pf}
If $C$ is a special $\lambda$-ring and $\Phi\colon A\to B$ a
$\lambda$-homomorphism, then for all
$x,y\in A$,
\begin{multline*}
\Phi(\lambda^n(xy)) = \lambda^n\Phi(xy) = \lambda^n(\Phi(x)\Phi(y)) \\
= P_n(\Phi(x),\ldots,\lambda^n\Phi(x),\Phi(y),\ldots,\lambda^n \Phi(y))\\
= P_n(\Phi(x),\ldots,\Phi(\lambda^n x),\Phi(y),\ldots, \Phi(\lambda^n y))\\
= \Phi(P_n(x,\ldots,\lambda^n x,y,\ldots,\lambda^n y)),
\end{multline*}
so
\begin{equation}\label{kerone}
\lambda^n(xy)-P_n(x,\ldots,\lambda^n x,y,\ldots,\lambda^n
y)\in\ker\Phi.\end{equation}
Similarly,
\begin{equation}\label{kertwo}
\lambda^m\lambda^n x - Q_{m,n}(x,\ldots,\lambda^{mn}
x)\in\ker\Phi.\end{equation}
Let $I$ denote the $\lambda$-ideal in $A$ generated by elements of type (\ref{kerone})
and (\ref{kertwo}).
Thus $I\subset \ker\Phi$.  Conversely, if $B=A/I$, (\ref{kerone}) and
(\ref{kertwo})
imply
 (\ref{ringhom}) and (\ref{lambdahom}) respectively.
Thus, the quotient map $A\to B$ is
universal.\end{pf}

We call the ring $B$, the {\bf specialization} of $A$.

Next we say a few words about the relationship between symmetric and
exterior powers.
In the case that $A$ is the Grothendieck group of vector spaces (resp.
projective modules,
resp. vector bundles), one-dimensional spaces (resp. invertible ideals,
resp. line bundles) correspond to elements $a$ such that $\lambda_t(a) = 1
+ at$.
If $b = a_1+\cdots+a_n$, where each $a_i$ satisfies this condition,
$$\lambda_t(b)=\prod_{i=1}^n\lambda_t(a_i)=\prod_{i=1}^n (1+a_i t).$$
Thus,
$$\prod_{i=1}^n (1+a_i t + a_i^2 t^2 + \cdots) = \prod_{i=1}^n \frac{1}
{1-a_i t}
= \lambda_{-t}(b)^{-1}.$$
We therefore define $\sigma^n(x)$ as the $n$th coefficient of
$\lambda_{-t}(x)^{-1}$
for all $x\in A$.  We note that if
\begin{equation}
\label{inversion}
\sigma_t(x)=\sum_{n=0}^\infty \sigma^n(x)t^n = \lambda_{-t}(x)^{-1},
\end{equation}
then
$$\sigma_t(x+y)=\lambda_{-t}(x+y)^{-1} =
\lambda_{-t}(x)^{-1}\lambda_{-t}(y)^{-1}
=\sigma_t(x)\sigma_t(y).$$
We conclude that
\begin{align*}
\sigma^0(x) &= 1\\
\sigma^1(x) &= x\\
\sigma^n(x+y) &= \sum_{i+j=n}\sigma^i(x)\sigma^j(y),
\end{align*}
so the $\sigma^n$ give a new $\lambda$-ring structure on $R$ which we
call the {\bf opposite} structure to $\{\lambda^n\}$.   

Several remarks on this construction are in order.  The automorphism on the
multiplicative group of power series with constant term $1$ given by
$f(t)\mapsto f(-t)^{-1}$ is an involution, so
the opposite of the opposite of a $\lambda$-structure is the structure itself.
By (\ref{inversion}), an element is virtually finite with respect to 
a $\lambda$-structure if and only if it is virtually finite with respect to the opposite
structure.
The opposite of a special $\lambda$-ring need not be special.  For example,
in $\bbZ$, which is special with respect to $\{\lambda^n\}$,
$\sigma^n(r) = \binom{r+n-1}{n}$, so the $\{\sigma^n\}$ counterpart of the identity
(\ref{ringhom}) for $n=2$, namely
$$\sigma^2(xy) = y^2 \sigma^2 x + x^2 \sigma^2 y - 2(\sigma^2 x)(\sigma^2 y)$$
does not hold.

The Newton polynomials, giving power sums in terms of elementary symmetric
functions,
allow us to define the {\bf Adams operations} $\Psi^n\colon A\to A$ on
any $\lambda$-ring.
In fact, we need these operations only for vector bundles over algebraic
varieties, where
they have been extensively studied.

\section{Special $\lambda$-ring $\overline{K}(X)$}

Let $X$ be a variety. It is well known that the usual $K$-theory (of algebraic
vector bundles) associates to $X$ a special $\lambda $-ring $K(X)$. We will need to
consider a different Grothendieck group $\overline{K}(X)$ of vector bundles, so that
the functor of global sections descends to a group homomorphism
$$H^0:\overline{K}(X)\to K[\Vect].$$
So let us take $\overline{K}(X)$ to be the abelian group generated by isomorphism
classes of algebraic vector bundles on $X$ with relations
$$[P]=[M]+[N],$$
whenever the vector bundles $M\oplus N$ and $P$ are isomorphic. Note that we do not
impose relations coming from general short exact sequences as in the usual $K$-theory. The
$\otimes $ operation makes $\overline{K}(X)$ a ring and the $\lambda$-operations are defined
in the usual way using the exterior powers:
$$\lambda ^i[P]=[\Lambda ^iP].$$

\begin{thm} The $\lambda$-ring $\overline{K}(X)$ is special.
\end{thm}

\begin{pf} The traditional way to prove that $K(X)$ is a special $\lambda$-ring uses the
splitting principle, which, in fact, is equivalent to the  $\lambda$-ring (in which every
element is of finite dimension) being special
\cite{FuLa}~Ch. 1. The usual method of splitting an algebraic vector bundle produces only a
short exact
sequence and therefore is not applicable to our group $\overline{K}(X)$. We will prove
the identities (\ref{ringhom}) and (\ref{lambdahom}) in $\overline{K}(X)$ 
by showing that for any $x,y\in \overline{K}(X)$
there exists a $\lambda$-homomorphism from a special $\lambda$-ring to $\overline{K}(X)$ such
that $x,y$ are contained in the image. For this we need a {\it free} special $\lambda$-ring
on two generators.

First recall the free special $\lambda $-ring on one generator. It has two standard
descriptions: as a ring of symmetric functions (in an infinite number of variables), and as
a direct sum of representation rings of the symmetric group. It is the second description
which is useful for our purposes, so we recall it.

Let $R_n$ be the representation ring of the symmetric group $S_n$,
with the convention that $S_0=\{ e\}$. That is, $R_n$ is a free
$\bbZ$-module with basis consisting of isomorphism classes of
irreducible (complex) representations of $S_n$. Put
$$R:=\bigoplus
_{n\geq 0}R_n.$$
 Thus $R$ has a $\bbZ$-basis consisting of pairs
$(n, \omega)$, where $n$ is a natural number and $\omega $ is an
irreducible representation of $S_n$. The ring structure is
uniquely determined by the requirement that 
$$(n_1,\omega_1)(n_2, \omega _2)=(n_1+n_2,\Ind_{S_{n_1}\times
S_{n_2}}^{S_{n_1+n_2}} \omega _1\otimes \omega _2).$$ 
The trivial $S_0=\{e\}$-module is the unit in the ring $R$. The $\lambda$-ring 
structure is defined as follows. Fix a basis element $(n, \omega )$ 
and a positive integer $r$. Consider the obvious action of the wreath product
$S_n^r\rtimes S_r$ on $\omega ^{\otimes r}$ and let $\Sign(S_r)$ be the sign representation
of $S_r$. Then $\omega ^{\otimes r}\otimes\Sign(S_r)$ is a left $S_n^r\rtimes S_r$-module.

Then
$$\lambda ^r(n,\omega )=(rn, \Ind _{S_n^r\rtimes S_r}^{S_{rn}}\omega ^{\otimes r}
\otimes \Sign(S_r)).$$
The special $\lambda $-ring $R$ is generated by the element $(1,\bbC)$: the elements
$\{ \lambda ^r(1, \bbC)\mid r\geq 1\}$ are algebraically independent in $R$.  The ring
$R$ is free in the following sense: given a special $\lambda$-ring $Q$ and an element
$x\in Q$ there exists a $\lambda $-homomorphism $f:R\to Q$ such that $f((1,\bbC))=x$.

Consider the special $\lambda $-ring $R^2:=R\otimes _{\bbZ}R$. It is the {\it free} special
$\lambda$-ring on two generators in the obvious sense. Naturally
$R^2$ has a $\bbZ$-basis consisting
of elements $((n_1,n_2), \omega _1\otimes \omega _2)$, where $\omega _i$ is an irreducible
representation of $S_{n_i}$, $i=1,2$; we regard $\omega_1\otimes\omega_2$
as a $S_{n_1}\times S_{n_2}$ representation in the usual way.
The $\lambda $-operations are similar to those in $R$:
$$\lambda ^r((n_1,n_2),\omega _1\otimes \omega _2 )=((rn_1,rn_2),
\Ind _{(S_{n_{1}}^r\times S_{n_{2}}^r)\rtimes S_r}^{S_{rn_1}\times S_{rn_2}}
\omega _1^{\otimes r}\otimes\omega _2^{\otimes r}\otimes \Sign(S_r)).$$

Let $M$, $N$ be (virtual) vector bundles on $X$.
Consider $M^{\otimes n_1}\otimes N^{\otimes n_2}$ as a (virtual) right 
$S_{n_1}\times S_{n_2}$-module.
We define the homomorphism
$\theta :R^2\longrightarrow \overline{K}(X)$ as follows.
$$\theta ((n_1,n_2), \omega _1\otimes \omega _2)=(M^{\otimes n_1}\otimes N^{\otimes n_2})
\otimes _{S_{n_1}\times S_{n_2}} (\omega _1\otimes \omega _2).$$
This is a ring homomorphism.
We only need to check that the $\lambda $-operations correspond under $\theta$; that is,
we need to compare
\begin{align*}
 \theta (\lambda ^r((n_1,n_2),&\omega _1\otimes \omega _2)) \\
 = 
 (M^{\otimes rn_1}\otimes &N^{\otimes rn_2})\otimes _{S_{rn_1}\times S_{rn_2}} \\
&( \bbC [S_{rn_1}\times S_{rn_2}]\otimes _{(S^r_{n_1}\times S_{n_2}^r)\rtimes S_r}
 (\omega _1^{\otimes r}\otimes \omega _2^{\otimes r}\otimes \Sign(S_r)))\\
 =  (M^{\otimes rn_1}\otimes &N^{\otimes rn_2})\otimes
 _{(S^r_{n_1}\times S_{n_2}^r)\rtimes S_r}
 (\omega _1^{\otimes r}\otimes \omega _2^{\otimes r}\otimes \Sign(S_r))\\
 \end{align*}
and
$$\begin{array}{rl}
 &  \Lambda ^r(\theta ((n_1,n_2),\omega _1\otimes \omega _2))\\
 = & \Lambda ^r((M^{\otimes n_1}\otimes N^{\otimes n_2})\otimes _{S_{n_1}\times S_{n_2}}
 (\omega _1\otimes \omega _2))\\
 = & ((M^{\otimes rn_1}\otimes N^{\otimes rn_2})\otimes _{S^r_{n_1}\times S^r_{n_2}}
 (\omega _1^{\otimes r}\otimes \omega _2^{\otimes r}))\otimes _{S_r}\Sign(S_r).\\
  \end{array}
$$
The following lemma implies the two are isomorphic.

\begin{lemma} Let $G$ be a group, $P$, $Q$ -- right and left $G$-modules respectively.
Consider the space $(P\otimes _G Q)^{\otimes r}=P^{\otimes r}\otimes _{G^r}Q^{\otimes r}$
as a right $S_r$-module.
There is
a canonical isomorphism of right $S_r$-modules
\begin{equation}
\label{genisom}
\alpha\colon P^{\otimes r}\otimes _{G^r}Q^{\otimes r}\stackrel{\sim}{\longrightarrow}
P^{\otimes r}\otimes _{G^r\rtimes S^r}(Q^{\otimes r}\otimes \bbC [S_r]).
\end{equation}
In particular there is a canonical isomorphism of vector spaces
\begin{equation}
\label{specisom}
\Lambda ^r((P\otimes _G Q)^{\otimes r})\stackrel{\sim}{\longrightarrow}
P^{\otimes r}\otimes _{G^r\rtimes S^r}(Q^{\otimes r}\otimes \Sign(S_r)).
\end{equation}
\end{lemma}

\noindent{\bf Proof of the lemma.}
Indeed,
$$\alpha\colon x\otimes y\mapsto x\otimes y\otimes 1,$$
and
$$p_1\otimes\cdots\otimes p_r\otimes q_1\otimes\cdots\otimes q_r\otimes \tau
\mapsto p_{\tau (1)}\otimes\cdots\otimes p_{\tau (r)}
\otimes q_{\tau (1)}\otimes\cdots\otimes q_{\tau (r)}$$
define mutually inverse $S_r$-linear maps.
This proves (\ref{genisom}).  Tensoring both sides over $S_r$ with $\Sign(S_r)$,
we get (\ref{specisom}) and therefore the proposition.
\end{pf}

So we have constructed a contravariant functor
$$\overline{K}:\{ \text{Varieties}\} \to \{ \text{Special $\lambda$-rings}\}.$$
The Adams operations $\Psi ^n$ are defined in $\overline{K}(X)$ in the usual way:
these are certain universal polynomials with integer coefficients in the operations
$\lambda ^i$. We have
$$\Psi ^n[\cL]=[\cL ^{\otimes n}],$$
if $\cL$ is a line bundle. Since the $\lambda $-ring $\overline{K}(X)$ is special,
the $\Psi ^n$ are $\lambda$-homomorphisms (see [AtiTa], Section 5).

Note that we have a well defined group homomorphism
$$H^0: \overline{K}(X)\to K_0[\Vect].$$

\medskip

\section{Motivic measures $\mu _n$}

Consider the ($\lambda $-) ring $\bbZ[M]$ as in Example~\ref{target}. We will freely consider
elements of the monoid $M$ either as polynomials with integer coefficients or as
(isomorphism classes of) graded vector spaces.

For a smooth connected projective variety $X$ of dimension $d$ define
$$\mu _n(X):=1+h^1_n(X)s+\cdots+h^d_n(X)s^d \in \bbZ[M],$$
where 
$$h^i_n(X)=\dim H^0(X,\Psi ^n\Omega ^i_X)\in \bbZ.$$

We constructed the measure $\mu =\mu _1$ in our previous paper [LaLu]. Explicitly
$$\mu (X)=1+h^{1,0}(X)s+\cdots+h^{d,0}(X)s^d$$
for a smooth projective irreducible $X$ of dimension $d$.

\begin{prop} For each $n\geq 1$ we have

i) $\mu _n(X)=\mu _n(\tilde {X})$ if $X$ and $\tilde{X}$ are birational,

ii) $\mu _n(X\times Y)=\mu _n(X)\mu _n (Y),$

iii) $\mu _n(\bbP ^k)=1$ for all $k\geq 0$.
\end{prop}

\begin{cor} The mapping $\mu _n$ extends (uniquely) to a ring homomorphism
$$\mu _n:K_0[\cV _{\bbC}]/(\bbA^1)\to \bbZ [M].$$
\end{cor}

\medskip

Indeed, this follows from [LaLu], Theorem 2.3.

\medskip

To prove the proposition, we need the following lemma:

\begin{lemma}
\label{boxtimes}
Let $X$ and $Y$ be varieties and
$E$ and $F$ vector bundles on $X$ and $Y$ respectively. Then
$$H^0(X\times Y,p^*E\otimes q^*F)=H^0(X,E)\otimes H^0(Y,F).$$
\end{lemma}

\begin{pf}
By the projection formula $p_*(p^*E\otimes q^*F)=E\otimes p_*q^*F$, where
$p_*q^*F$ is the trivial bundle on $X$ with fiber $H^0(Y,F)$. Hence
\begin{align*}
H^0 (X\times Y,p^*E\otimes q^*F) &= H^0(X, p_*(p^*E\otimes q^*F))\\
 &= H^0(X, E\otimes p_*q^*F) = H^0(X,E)\otimes H^0(Y,F).
\end{align*}
\end{pf}
\smallskip
\goodbreak
\noindent{\bf Proof of proposition.}

i) Since $\Psi ^n$ is a polynomial in operations $\lambda ^j$, it suffices to prove
that $H^0(X,(\Omega _X^i)^{\otimes n}))=H^0(\tilde{X}, (\Omega _{\tilde{X}}^i)^{\otimes n})$.
But this is well known (the proof is the same as that of \cite{Hart}~II~8.19).

ii) Consider the projections $X\stackrel{p}{\leftarrow}X\times Y\stackrel{q}{\rightarrow} Y$.
We have
$$\Omega ^1_{X\times Y}=p^*\Omega _X^1\oplus q^*\Omega _Y^1.$$
Hence
$$\Omega ^m_{X\times Y}=\bigoplus_{i+j=m}p^*\Omega ^i_X\otimes q^*\Omega _Y^j.$$
As $\bar K(X)$ is a special $\lambda$-ring,
\begin{equation*}
\begin{split}
\Psi^n\Omega _{X\times Y}^m & =  \bigoplus _{i+j=m}\Psi ^np^*\Omega _X^i\otimes
\Psi ^nq^*\Omega _Y^j\\
& =  \bigoplus _{i+j=m}p^*\Psi ^n\Omega _X^i\otimes
q^*\Psi ^n\Omega _Y^j.\\
\end{split}
\end{equation*}
It follows that
$$H^0(X\times Y,\Psi ^n\Omega _{X\times Y}^m)=\bigoplus_{i+j=m}H^0(X, \Psi ^n\Omega _X^i)
\otimes H^0(Y, \Psi ^n\Omega _Y^j),$$
i.e. $\mu _n(X\times Y)=\mu _n(X)\mu_n(Y)$.

iii) Since $\bbP ^k$ is  birational to $(\bbP ^1)^k$, by i),ii) above it suffices to prove
that $\mu _n(\bbP ^1)=1$. We have
$$\Psi ^n\Omega ^i_{\bbP ^1}=
\begin{cases}
\cO _{\bbP ^1},\hbox{ if $i=0$,}\\
\cO _{\bbP ^1}(-2n),\hbox{ if $i=1$,}\\
0,\hbox{ if $i>1$.}\\
\end{cases}
$$
and therefore $\mu _n(\bbP ^1)=1\in M$. This proves the proposition.
\qed\smallskip

\begin{prop} Let $X$ be a variety of dimension $\leq 2$. Then
$$\mu (\Sym ^mX)=\lambda ^m(\mu (X)).$$
\end{prop}

\begin{pf} The class of $X$ in $K_0[\cV _{\bbC}]$ is
a linear combination of classes of smooth projective varieties of dimension 0,1 and 2.
Thus we may assume that $X$ is smooth projective. If $\dim X=0$ the assertion is
trivial. We will prove the assertion in case $\dim X=2$ (the case of curves is similar
and simpler since $\Sym ^nX$ is smooth if $\dim X=1$).
So let $X$ be a smooth projective surface.

We will use the smooth variety $\Hilb ^mX$---the Hilbert scheme of zero-dimensional
subschemes of $X$.  It follows immediately from
the main theorem of L.~Goettsche \cite{Go} that the classes of  $\Hilb ^mX$ and $\Sym ^mX$
 are equal in $K_0[\cV _{\bbC}]/[\bbA ^1]$.  In particular,
 $$\mu _n (\Hilb ^mX)=\mu _n(\Sym ^mX).$$
By Lemma~\ref{omegai} below 
$H^0(\Hilb ^mX,\Omega_{\Hilb^m X}^i)=H^0(X^m,\Omega_{X^m}^i)^{S_m}$. Thus we need to
prove the following equality in $M$:
\begin{multline*}
1+H^0(X^m,\Omega ^1)^{S_m}s+\cdots+H^0(X^m,\Omega ^{2m})^{S_m}s^{2m}\\
=\lambda ^m(1+H^0(X,\Omega ^1)s+ H^0(X,\Omega ^2)s^2).\\
\end{multline*}

Recall from Example~\ref{vvs} 
that
$$\lambda ^j(V s^i)=\begin{cases}
\Sym ^j(V)s^{ij},\ \text{if $i$ is even},\\
\Lambda ^j(V)s^{ij},\ \text{if $i$ is odd}.\\
\end{cases}
$$ 
Let $\sum_iV_is^i\in M$ be a graded vector space. Then
$$\lambda ^m\bigl(\sum_iV_is^i\bigr)=\sum_i\sum
_{j_1n_1+\cdots+j_kn_k=i}\lambda ^{n_1}(V_{j_1}) \otimes \cdots\otimes
\lambda ^{n_k}(V_{j_k})s^i.$$ 
Let us prove the following general
lemma.  Let $Y$ be a smooth projective variety and $F$ a vector
bundle on $Y$ of rank $r$. Fix $m\geq 1$ and denote by $p_i:Y^m\to
Y$ the projection to the $i$th factor. Put $F_i:=p_i^*F$, $\cF:= F_1\oplus\cdots\oplus F_m$.

\medskip

\begin{lemma} There is a natural isomorphism of graded vector spaces
$$\sum_jH^0(Y^m,\Lambda ^j(\cF))^{S_m}s^j
=\lambda ^m\bigl(\sum_jH^0(Y,\Lambda ^jF)s^j\bigr).$$
\end{lemma}

\begin{pf}
We have
$$\Lambda ^\bullet (\cF)=\bigoplus_{(j_1,\ldots,j_m)}\Lambda ^{j_1}F_1\otimes \cdots\otimes
\Lambda ^{j_m}F_m,$$
%
By Lemma~\ref{boxtimes},
$$H^0(Y^m, \Lambda ^{j_1}F_1\otimes\cdots\otimes
\Lambda ^{j_m}F_m)\simeq \bigotimes_{i=1}^m H^0(Y^m,\Lambda ^{j_i}F_i)\simeq
\bigotimes_{i=1}^m H^0(Y,\Lambda^{j_i}F).$$

Put $\cF _{j_1\cdots j_m}=\Lambda ^{j_1}F_1\otimes \cdots\otimes \Lambda ^{j_m}F_m$.
Consider the $S_m$-action on the sheaf $\Lambda ^\bullet(\cF)$. The sheaves
$\cF _{j_1\cdots j_m}$ and $\cF _{j^\prime _1\cdots j^\prime _m}$ are in the same $S_m$-orbit
if and only if the multisets  $\{j_1,\ldots,j_m\}$ and $\{j^\prime_1,\ldots,j^\prime_m\}$ are equal.
For each multiset $\{j_1,\ldots,j_m\}$ choose a representative $(j_1,\ldots j_m)$ of the
corresponding $S_m$-orbit. Put $G _{j_1\cdots j_m}=\Stab_{S_m}(j_1,\ldots j_m)$.
Then
$$W=H^0(Y^m,\Lambda ^\bullet \cF)^{S_m}
=\sum_{\{j_1,\ldots,j_m\}}H^0(Y^m,\cF _{j_1\cdots j_m})^{G_{j_1\cdots j_m}}.$$

Consider the space
$$H^0(Y^m,\cF _{j_1\cdots j_m})^{G_{j_1\cdots j_m}}=
(\otimes _sH^0(Y^m,\Lambda^{j_s}F_s))^{G_{j_1\cdots j_m}}.$$
Assume for simplicity of notation that the multiset $\{j_1,\ldots,j_m\}$ contains $k$ different
elements:
$$j_1=\cdots=j_{t_1}\neq j_{t_1+1}=\cdots=j_{t_1+t_2}\neq \cdots=j_{t_1+\cdots+t_k}.$$
Then $G_{j_1\cdots j_m}=S_{t_1}\times \cdots \times S_{t_k}$. For example, the subgroup $S_{t_1}$
acts by permuting factors in 
$\Lambda ^{j_1}F_1\otimes \cdots\otimes \Lambda ^{j_1}F_{t_1}$.
A transposition of two factors in this tensor product corresponds to a product of
$j_1$ transpositions on the level of $\cF^{\otimes j_1 t_1}$.  Therefore,
$$H^0(Y^m, \Lambda ^{j_1}F_1\otimes \cdots\otimes \Lambda ^{j_{t_1}}F_{t_1})^{S_{t_1}}
\simeq
\begin{cases}
\Sym ^{t_1}H^0(Y,\Lambda ^{j_1}F),\hbox{ if $j_1$ is even,}\\
\Lambda ^{t_1}H^0(Y,\Lambda ^{j_1}F),\hbox{ if $j_1$ is odd,}\\
\end{cases}
$$
The space $H^0(Y^m,\cF _{j_1\cdots j_m})^{G_{j_1\cdots j_m}}$ is the tensor product of $k$ factor
of the form
$\Sym ^{t_p}H^0(Y,\Lambda ^{l_p}F)$ (if $l_p$ is even) or $\Lambda ^{t_p}H^0(Y,\Lambda ^{l_p}F)$
(if $l_p$ is odd). Also the degree of this space is equal to  $\sum_s j_s$.
This proves the lemma.
\end{pf}

Now apply the lemma with $Y=X$ and $F=\Omega_X^1$ to get the proposition.
\end{pf}

\medskip

\begin{remark} Unfortunately, the assertion of the last proposition is no longer true if
we replace the measure $\mu $ by $\mu _n$ for $n\geq2$. A counter example is provided by a
smooth projective curve of genus $\geq 2$. Indeed, then $\Sym ^mX$ is stably birational
to the Jacobian of $X$ and hence $\mu _n(\Sym ^mX)=\mu (\Sym ^mX)$. On the other hand
it is clear that $\mu _n(X)$ and hence $\lambda ^m(\mu _n(X))$ depends on $n$.
\end{remark}

\medskip

\begin{question} Is $\mu :K_0[\cV _{\bbC}]\to \bbZ[M]$ is a $\lambda$-homomorphism?
It seems likely to us that the answer is affirmative.
We could prove it if we knew that for any smooth projective
variety $Z$ and any $m\geq 1$ the class of $\Sym^mZ$ is equal in $K_0[\cV _{\bbC}]/[\bbA ^1]$
to the class of a resolution of $\Sym ^mZ$. 
So one might generalize and ask the following question.  Let $X$ be a 
non-singular complex projective variety with an action
by a finite group $G$, and $Y$ a non-singular projective variety birationally
equivalent to $X/G$.  Is it always true that
$$[X/G]\equiv [Y]\pmod{[\bbA ^1]}?$$

\end{question}

\medskip

\section{Irrationality theorems}

Let $X$ be a smooth projective surface. Let $m\geq 1$ be an integer and $\Hilb ^mX$
be the Hilbert scheme of zero-dimensional subschemes of $X$ of length $m$. It is
well known that the natural (Hilbert-Chow) 
morphism $g: \Hilb ^mX\to \Sym ^mX$ is a resolution of
singularities which is a semismall map ([Na]~1.15,\,6.10). It follows immediately from
the main theorem of L.~Gottsche \cite{Go} that the classes of  $\Hilb ^mX$ and $\Sym ^mX$
 are equal in $K_0[\cV _{\bbC}]/[\bbA ^1]$. Thus for all $n\geq 1$
 $$\mu _n (\Hilb ^mX)=\mu _n(\Sym ^mX).$$

Let us introduce some notation.  Consider the pullback diagram
$$\begin{array}{ccc}
Y & \stackrel{f}{\longrightarrow} & X^m \\
\pi \downarrow & & \downarrow \pi \\
\Hilb ^mX & \stackrel{g}{\longrightarrow} & \Sym ^mX,
\end{array}
$$
where the vertical arrows are quotient morphisms 
by the $S_m$-action. Let $\Sym ^mX_*$ denote
the open subspace of $\Sym ^mX$ consisting of 
all multisets $\{x_1,\ldots,x_m\}$ in which at least
$m-1$ points are distinct. Pulling back the previous diagram to $\Sym ^mX_*$,
we have
$$\begin{array}{ccc}
Y_* & \stackrel{f}{\longrightarrow} & X^m_* \\
\pi \downarrow & & \downarrow \pi \\
\Hilb ^mX_* & \stackrel{g}{\longrightarrow} & \Sym ^mX_*.
\end{array}
$$

Denote by $\Delta _*$ the intersection of the ``fat'' diagonal $\Delta \subset X^m$ with
$X^m_*$. Then $Y_*$ is the blowup of $X^m_*$ along the smooth subvariety $\Delta _*$ which
is of codimension 2. So $Y_*$ is smooth.
 Let $E\subset Y_*$ be the exceptional divisor. The map $\pi :Y_*
\to \Hilb ^mX_*$ is unramified away from $E$ and has ramification of degree 2 along $E$.
Since the map $g$ is semismall, the complement of the
open subset $\Hilb ^mX_*\subset \Hilb ^mX$
has codimension $\geq 2$. Hence also $Y\setminus Y_*$ has codimension $\geq 2$ in $Y$.

Finally, let $\Sym^m _+X$ be the image of $X^m - \Delta $.  We pull back the previous diagram to
$\Sym ^m_+X$:
\begin{equation}
\label{pbdiag}
\begin{array}{ccc}
Y_+ & \stackrel{f}{\longrightarrow} & X^m_+ \\
\pi \downarrow & & \downarrow \pi \\
\Hilb ^mX_+ & \stackrel{g}{\longrightarrow} & \Sym ^mX_+.
\end{array}
\end{equation}
The vertical maps $\pi$ are \'etale and the horizontal maps $f,g$ are isomorphisms.
The complement of the open subset $X^m_+\subset X^m$ has codimension 2.

\medskip

\begin{lemma} 
\label{multione}
For any $m,n\geq 1$ we have
$$H^0(\Hilb ^mX,(\Omega^1_{\Hilb ^mX})^{\otimes n})
\subset H^0(X^m,(\Omega^1_{X^m})^{\otimes n})^{S_m}.$$
\end{lemma}

\begin{pf} Consider the diagram (\ref{pbdiag}). We have
\begin{multline*}H^0(\Hilb ^mX,(\Omega^1_{\Hilb ^mX})^{\otimes n})\subset H^0(\Hilb ^mX_+,(\Omega^1_{\Hilb ^mX_+})^{\otimes n})\\
= H^0(Y_+,(\Omega^1_{Y_+})^{\otimes n})^{S_m}
= H^0(X_+^m,(\Omega^1_{X_+^m})^{\otimes n})^{S_m}.\\
\end{multline*}
Moreover, 
$$H^0(X_+^m,(\Omega^1_{X_+^m})^{\otimes n})^{S_m}
= H^0(X^m,(\Omega^1_{X^m})^{\otimes n})^{S_m},$$
since $X^m_+$ has complement of codimension 2 in $X^m$.
\end{pf}

\medskip

\begin{lemma} 
\label{plurig}
For any $m,n\geq 1$ we have
$$H^0(\Hilb ^mX,\omega_{\Hilb ^mX} ^{\otimes n})
 = H^0(X^m, \omega_{X^m} ^{\otimes n})^{S_m}.$$
\end{lemma}

\begin{pf} It suffices to prove that
$$H^0(\Hilb ^mX_*,\omega_{\Hilb ^mX_*}^{\otimes n})
 = H^0(X^m_*, \omega_{X^m_*}^{\otimes n})^{S_m}.$$

We have the natural injective maps
$$\alpha :\pi ^*\omega _{\Hilb ^mX_*}\to \omega _{Y_*},
\quad \beta:f^* \omega _{X^m_*}\to \omega _{Y_*}.$$
It suffices to prove that $im(\alpha )=im(\beta)$.  Indeed, 
$$H^0(Y_*,\pi ^*\omega _{\Hilb ^mX_*})^{S_m}=H^0(\Hilb ^mX_*,\omega _{\Hilb ^mX_*}).$$
Both $\alpha$ and $\beta$ are surjective away from $E$. So it remains to analyze the
maps $g$ and $\pi$ near $E$.

Choose a point $p=(a_1,a_2=a_1,a_3,\ldots,a_m)\in \Delta _*$ and $q\in g^{-1}(p)$. There exist
local (analytic) coordinates $x_1,x_2,\ldots,x_{2m}$ near $p$ and local coordinates
$y_1,\ldots,y_{2m}$ near $q$ so that
$$g^*(x_1)=y_1,\ g^*(x_2)=y_1y_2, \ g^*(x_3)=y_3,\  \ldots \ g^*(x_{2m})=y_{2m}.$$
Then $y_1=0$ is the local equation of $E$.
Thus 
$$g^*(dx_1\wedge \cdots\wedge dx_{2m})=y_1dy_1\wedge \cdots\wedge y_{2m}$$
and
$$\beta ((dx_1\wedge \cdots\wedge dx_{2m})^{\otimes n})=y_1^n(dy_1\wedge \cdots\wedge y_{2m})^{\otimes n},$$
so that $im(\beta)=\omega _{Y_*}^{\otimes n}(-nE)$.
Similarly, we can define local 
coordinates $z_1,\ldots,z_{2m}$ near $\pi (q)\in \Hilb ^mX_*$, so that
$$\pi ^*(z_1)=y_1^2,\ \pi ^*(z_2)=y_2,\ \ldots\ \pi ^*(z_{2m})=y_{2m}.$$
Therefore 
$$g^*(dz_1\wedge \cdots \wedge dz_{2m})=2y_1dy_1\wedge \cdots \wedge dy_{2m}$$
and
$$\alpha ((dz_1\wedge \cdots\wedge dz_{2m})^{\otimes n})=2^ny_1^n
(dy_1\wedge \cdots\wedge y_{2m})^{\otimes n}.$$
That means $im(\alpha)=\omega _{Y_*}(-nE).$

\end{pf}

\medskip

\begin{lemma} 
\label{symm}
For any $m,n\geq 1$ we have
$$H^0(X^m,\omega_{X^m}^{\otimes n})^{S_m}=\Sym ^m H^0(X,\omega_X^{\otimes n}).$$
\end{lemma}

\begin{pf}
Let $p_i:X^m\to X$ denote the projection on the $i$th factor. Then
$$\omega ^{\otimes n} _{X^m}
=(p_1^*\omega _X\otimes\cdots\otimes p_m^*\omega _X)^{\otimes n}\simeq
p_1^*\omega ^{\otimes n}_X\otimes\cdots\otimes p_m^*\omega ^{\otimes n}_X,$$
and hence
$$H^0(X^m, \omega_{X^m}^{\otimes n})=H^0(X,\omega_X^{\otimes n})^{\otimes m}.$$
The $S_m$-action permutes the factors and
$$H^0(X^m,\omega_{X^m}^{\otimes n})^{S_m}=\Sym^m H^0(X,\omega_X^{\otimes n}),$$
since $\dim X$ is even.
\end{pf}

\medskip

\begin{lemma} 
\label{omegai}
For any $m,i\geq 1$ we have  
$$H^0(\Hilb ^mX,\Omega_{\Hilb ^mX}^i) = H^0(X^m,\Omega_{X^m}^i)^{S_m}.$$
\end{lemma}

\begin{pf} Since $Y_*$ is a blowup of $X^m_*$ along a smooth subvariety we have
$H^0(Y_*,\Omega_{Y^*}^i)=H^0(X^m_*,\Omega_{X^m}^i)$. 
Also $H^0(X^m_*,\Omega_{X^m_*}^i)=H^0(X^m,\Omega_{X^m}^i)$,
since $X^m_*$ has complement of codimension $\geq 2$ in $X^m$. Thus we must show that
$$H^0(\Hilb ^mX,\Omega_{\Hilb^m X}^i)= H^0(Y_*,\Omega_{Y_*}^i)^{S_m}.$$
Let $q\in Y_*$, $G_q=\Stab(q)\subset S_m$. It suffices to show that
$$(\pi ^{-1}\Omega _{\Hilb ^mX_*})_q=(\Omega _{Y_*,q})^{G_q}.$$
If $q\in Y_+$, then $G_q=1$ and $(\pi ^{-1}\Omega _{\Hilb ^mX_*})_q=(\Omega _{Y_*,q}).$
Suppose $q\in E$ and hence $G_q=\bbZ /2\bbZ$.
Consider local coordinates $y_1,\ldots,y_{2m}$ at $q$ and local
coordinates $z_1,\ldots,z_{2m}$ at $\pi (q)$ as in the proof of Lemma~\ref{plurig} above, so that
$$\pi ^*(z_1)=y_1^2,\ \pi ^*(z_2)=y_2,\ \ldots\ \pi ^*(z_{2m})=y_{2m}.$$
An element $y_1^kdy_1$ is $G_q$-invariant if and only if $k$ is odd.
The lemma follows.
\end{pf}

\begin{prop} 
\label{absval}
Fix integers $n,\ i\geq 1$.
Then the absolute value of the integer $h_n^i(\Hilb^mX)$ is bounded independently of $m$.
\end{prop}

\begin{pf} Recall that $h_n^i(Z)=H^0(Z,\Psi ^n\Omega ^i_Z)$. The Adams operation $\Psi ^n$
is a polynomial in the operations $\lambda ^i$. Hence it suffices to prove that
for all $n\ge 0$,
the dimension of the space $H^0(\Hilb ^mX, (\Omega^1_{\Hilb^m X})^{\otimes n})$
is bounded independently of $m$. By Lemma~\ref{multione} above
$$H^0(\Hilb ^mX, (\Omega^1_{\Hilb^m X} )^{\otimes n})
\subset H^0(X^m, (\Omega^1_{X^m} )^{\otimes n})^{S_m}.$$

We will prove that $\dim H^0(X^m, (\Omega^1_{X^m} )^{\otimes n})^{S_m}$ is bounded independently of $m$.

Let $p_i:X^m\to X$ be the projection to the $i$th factor. Denote $F_i=p_i^*\Omega^1_X$.
Then $\Omega^1_{X^m}=\bigoplus_{i=1}^m F_i$ and
$$(\Omega^1_{X^m})^{\otimes n}
=\bigoplus_{i_1,\ldots,i_n\in \{1,\ldots,m\}}F_{i_1}\otimes \cdots\otimes F_{i_n}.$$
The $S_m$-action permutes the summands $F_{i_1}\otimes \cdots\otimes F_{i_n}$ and the orbits
correspond to partitions $P$ of the set $\{1,\ldots,n\}$. Thus
$$H^0(X^m,\Omega ^{\otimes n})^{S_m}=\bigoplus _PH^0(X^m,F_P)^{\Stab(P)},$$
Where $F_P$ is one of the summands in the orbit, corresponding to $P$. Fix a partition $P$.
Assume for simplicity of notation that $P$ divides $\{1,\ldots,n\}$ into $k$ segments $P_1,\ldots,
P_k$, where $a<b$ for each $a\in P_s$, $b\in P_t$ if $s<t$. Let $\alpha _i=|P_i|$. Then
$$H^0(X^m,F_P)=H^0(X,\Omega ^{\otimes \alpha _1})\otimes \cdots\otimes
H^0(X,\Omega ^{\otimes \alpha _k}).$$
Therefore $\dim H^0(X^m,F_P)$ is bounded independently of $m$.
\end{pf}

\medskip

\begin{thm} 
\label{negative}
Let $X$ be a smooth projective surface of Kodaira dimension $\geq 0$. Then
the zeta-function $\zeta _X(t)$ is not pointwise rational.
\end{thm}

\begin{pf} Let $F$ denote the field of fractions of $\bbZ[M]$. We will show that there
exists $n\geq 1$ such that the power series
$$1+\sum_{m=1}^{\infty}\mu _n(\Sym^mX)t^m\in F[[t]]$$
is not rational.

First we claim that there exists $n\geq 1$ such that 
$$h^{2d}_n(\Hilb ^mX)=\dim H^0(\Hilb ^mX, \omega_{\Hilb^m X} ^{\otimes n})>0$$
for all
$m\geq 1$. Indeed, by our assumption on $X$ there exists $n\geq 1$,
such $H^0(X,\omega ^{\otimes n})\neq 0$ and by Lemmas \ref{plurig} and \ref{symm},
$$H^0(\Hilb ^mX, \omega_{\Hilb^m X}^{\otimes n})=\Sym ^mH^0(X,\omega_X^{\otimes n})$$
for all $m\geq 1$.
Fix one such $n$ and consider the
motivic measure $\mu _n: K_0[\cV _{\bbC}]\to \bbZ[M]\subset F$. We have
$\mu _n (\Sym ^mX)=\mu _n(\Hilb ^mX)$. It follows that $\mu _n(\Sym ^mX)\in M$ is a
polynomial with constant term 1 and leading term $h^{2d}_n(\Hilb ^mX)s^{2d}$.

Consider the group completion $G$ of the monoid $M$. By Lemma~\ref{monoidm},
$G$ is a free abelian group.
Hence we may apply Proposition~\ref{monomial} to the power series
$1+\sum_{m=1}^{\infty}\mu _n(\Sym^mX)t^m$.
Denote $g_m=\mu _n(\Sym ^mX)\in M$.
Assume that this power series is rational. Then there exist
$k,i\geq 1$ and an element $g\in G$ such that
$$g_{i+(\alpha +1)k}=gg_{i+\alpha k}$$
for all $\alpha \geq 0$.  The element $g$ is a rational function in $s$. Since all
coefficients $g_m$ are nonzero polynomials in $s$ 
it follows that $g$ is also a nonzero polynomial.
Note that $g$ in not a monomial since the degrees of the polynomials $g_m$ grow and their
constant term is 1. But then the coefficients of fixed powers of $s$ in the $g_m$ cannot
stay bounded, which contradicts Proposition~\ref{absval}.
\end{pf}

\medskip

\section{The special Grothendieck $\lambda$-ring of varieties.}

As pointed out in Example~\ref{source}, the symmetric power operations $\Sym^n$ define a
$\lambda$-ring structure on the Grothendieck ring $\groth{\bbC}$.
The zeta-function $\zeta_X(t)$ looks formally like a universal
$\lambda$-homomorphism,
but in fact it is not a ring homomorphism at all.  To see this, it suffices to
note that $\zeta_X(t)$ is a ratio of polynomials with constant term 1 for every curve $X$;
thus (\ref{prodform}) implies that $\zeta_X(t)\cdot\zeta_Y(t)$ is rational, while
Theorem~\ref{negative} implies $\zeta_{X\times Y}(t)$ is irrational whenever
$X$ and $Y$ both have positive genus.

Let $\{\lambda^n\}$ denote the $\lambda$-structure opposite to $\{\Sym^n\}$.  
This seems to be the natural choice of $\lambda$-structure on $\groth{\bbC}$
insofar as $X\mapsto \Sym^n X$ behaves like a symmetric power map, for instance
on cohomology.  The choice does not affect which classes are virtually finite, but it
will make a difference when we specialize, since, as we have seen, specialization 
does not commute with taking opposites.

\begin{lemma}
If $X$ is a variety over $\bbC$ whose image in $\groth{\bbC}$ is
virtually finite, then $\zeta_X(t)$ is globally rational.
\end{lemma}

\begin{pf}
Let $[X] = y - z$, where $y$ and $z$ are finite.  As $\lambda_t$ is a ring
homomorphism,
$$\zeta_X(T) = \lambda_{-t}([X])^{-1} = \lambda_{-t}(y)^{-1}\lambda_{-t}(z)$$
is globally rational in $\groth{\bbC}$.
\end{pf}

We have the following variant of Theorem~\ref{onedim}:

\begin{prop}
The class of any $1$-dimensional variety in $\groth{\bbC}$ is virtually
finite.
\end{prop}

\begin{pf}
If  $Y$ is a closed subvariety of $X$ with complement
$U$, then
if any two of $X$, $Y$, and $U$ are virtually finite, the third is so as
well.
A point $Z$ is finite dimensional (in fact $1$-dimensional in the $\lambda$-ring sense).
Thus, we reduce exactly as in the proof of Theorem~\ref{onedim}
to the case of a single projective non-singular curve $X$.
Now,
$$\lambda_t([{\mathbb P}^1]) = (1+t)(1+[\A^1]t),$$
so ${\mathbb P}^1$ is $2$-dimensional.  Writing $X = {\mathbb P}^1 - ({\mathbb P}^1-X)$,
(\ref{curveNumer}) implies
$$\lambda_t({\mathbb P}^1-X) =
(1+t)(1+\A^1t)\lambda_t(X)^{-1} = (1+t)(1+\A^1t)\zeta_X(-t) \in\groth{K}[t].$$
The proposition follows.\end{pf}

\begin{cor}
\label{vfbirat}
The virtual finiteness of a complex surface $X$ depends only on
the birational class of $X$.
\end{cor}

\begin{defn}
We call the specialization of the Grothendieck ring $\groth{K}$ with
respect to $\lambda^n$
the {\bf Grothendieck $\lambda$-ring} of $K$ and denote it
$\grsig{K}$.
The image of $\zeta_X(t)$ in $\grsig{K}$ will be denoted $\ze X$
\end{defn}

In any special $\lambda$-ring, the set of virtually finite elements is
clearly a
$\lambda$-subring.  In particular, we have:

\begin{prop}
For every $X/\bbC$ which is virtually finite in
$\groth{\bbC}$
and every positive integer $n$, $\Sym^n X$ is again virtually finite in
$\groth{\bbC}$.
\end{prop}

\begin{pf}
The identities relating symmetric and exterior powers show that $\Sym^n X$
lies in the
$\lambda$-subring of $\grsig{\bbC}$ generated by $X$.
\end{pf}

\begin{prop}
Every principally polarized abelian surface $X/\bbC$ is virtually finite in
$\grsig{\bbC}$.
\end{prop}

\begin{pf}
It is well known that every principally polarized abelian variety of
dimension $\leq 3$
is a product of Jacobian varieties.  (We do not know who first made this
observation, but
it is an immediate consequence of Torelli's theorem and Baily's theorem
\cite{Bai}.)
It suffices, therefore, to prove that the Jacobian $J$ of a genus $2$
curve $X$ is virtually
finite.   The map $\Sym^2 X\to J$ is a birational equivalence, so the
proposition follows from
Corollary~\ref{vfbirat}.
\end{pf}

We remark that
the fact that sufficiently high symmetric powers of a non-singular
projective curve
are projective space bundles over its Jacobian variety does not
immediately imply
virtual finiteness of Jacobians, since it is not obvious that the virtual
finiteness
of a projective space bundle or even a vector bundle over a given variety
implies the
virtual finiteness of that variety.  However, it is easy to prove the
pointwise rationality
of $\ze J$ for Jacobians $J$.

\begin{question}
Is $\ze A$ rational (globally or pointwise) for all abelian varieties $A$?
\end{question}

More optimistically, we might ask:

\begin{question}
Is $\ze X$ rational for all varieties $X$?  Is $\grsig{\bbC}$ finite dimensional in
the sense of Definition~\ref{fdim}?

\end{question}

\end{document}